\documentstyle[11pt,psfig,amssymb,amsmath]{article}
\newcommand{\be}{\begin{equation}} 
\newcommand{\ee}{\end{equation}}   
\newcommand{\bea}{\begin{eqnarray}}
\newcommand{\eea}{\end{eqnarray}}

\newtheorem{theorem}{Theorem}

\newtheorem{lemma}{Lemma}

\newtheorem{definition}{Definition}   

\def\1#1{^{(#1)}}

\begin{document}
%\DeclareGraphicsExtensions{.jpg,.pdf,.png}
\title{Analytic Representations in the 3-dim Frobenius Problem}

\author{Leonid G. Fel\\
\\Department of Civil and Environmental Engineering,\\
Technion, Haifa 3200, Israel\\
\vspace{-.3cm}  
\\{\sl e-mail: lfel@techunix.technion.ac.il}}
\date{\today}
\maketitle
\def\be{\begin{equation}}
\def\ee{\end{equation}}
\def\p{\prime}  
\begin{abstract}
We consider the Diophantine problem of Frobenius for semigroup ${\sf S}\left(
{\bf d}^3\right)$ where ${\bf d}^3$ denotes the tuple $(d_1,d_2,d_3)$, $\gcd(
d_1,d_2,d_3)=1$. Based on the Hadamard product of analytic functions 
\cite{horn90} we have found the analytic representation for the diagonal 
elements $a_{kk}\left({\bf d}^3\right)$ of the Johnson's matrix of minimal 
relations \cite{john60} in terms of $d_1,d_2,d_3$. Bearing in mind the results 
of the recent paper \cite{fel04} this gives the analytic representation for 
the Frobenius number $F\left({\bf d}^3\right)$, genus $G\left({\bf d}^3\right)$
and the Hilbert series $H({\bf d}^3;z)$ for the semigroups ${\sf S}\left({\bf 
d}^3\right)$. This representation does complement the Curtis' theorem 
\cite{curt90} on the non--algebraic representation of the Frobenius number 
$F\left({\bf d}^3\right)$. We also give a procedure to calculate the diagonal 
and off--diagonal elements of the Johnson's matrix.
\vskip .3cm
\noindent
\begin{tabbing}
{\sf Key words}:\hspace{.1in}  \=  Semigroups, Frobenius problem, 
    Hilbert series of a graded ring,\\
\> Hadamard product of analytic functions.
\end{tabbing}
\end{abstract}

%\tableofcontents
\newpage
\section{Introduction}\label{sect1}
Let ${\sf S}\left(d_1,\ldots,d_m\right)\subset {\mathbb N}$ be the numerical 
semigroup generated by a minimal set 
\footnote{The set $\{d_1,\ldots,d_m\}$ is called {\em minimal} if there are no 
nonnegative integers $b_{i,j}$ for which the following linear dependence holds
$d_i=\sum_{j\neq i}^mb_{i,j}d_j,\;\;b_{i,j}\in\{0,1,\ldots\}$ for any $i\leq m$.
It is classically known \cite{baru97} that $d_1\geq m$.}
of integers $\{d_1,\ldots,d_m\}$ such that
\begin{equation}
m\leq 
d_1<\ldots<d_m\;,\;\;\;\gcd(d_1,\ldots,d_m)=1\;,\;\;\;m\geq 2\;.\label{definn0}
\end{equation}
For short we denote the tuple $(d_1,\ldots,d_m)$ by ${\bf d}^m$. The least 
positive integer ($d_1$) belonging to ${\sf S}\left({\bf d}^m\right)$ is called 
{\em the multiplicity}. The number ($m$) of the minimal generators of ${\sf S}
\left(d_1,\ldots,d_m\right)$ is called {\em the embedding dimension}. {\em The 
conductor} $c\left({\bf d}^m\right)$ of ${\sf S}\left({\bf d}^m\right)$ is 
defined by $c\left({\bf d}^m\right):=\min\left\{s\in {\sf S}\left({\bf d}^m
\right)\;|\;s+{\mathbb N}\cup\{0\}\subset {\sf S}\left({\bf d}^m\right)\right\}
$. {\em The genus} $G\left({\bf d}^m\right)$ of ${\sf S}\left({\bf d}^m\right)$ 
is defined as the cardinality ($\#$) of its complement $\Delta$ in ${\mathbb N}
$, i.e. $\Delta\left({\bf d}^m\right)={\mathbb N}\setminus {\sf S}\left({\bf d}^
m\right)$ and
\begin{eqnarray}
G\left({\bf d}^m\right):=\#\Delta\left({\bf d}^m\right)\;.\label{defgen}
\end{eqnarray}
Introduce the generating function $\Phi\left({\bf d}^m;z\right)$ for the set
$\Delta({\bf d}^m)$ of unrepresentable integers
\begin{eqnarray}
\Phi\left({\bf d}^m;z\right)=\sum_{s\;\in\;\Delta\left({\bf d}^m\right)}z^s\;.
\label{def1}
\end{eqnarray}
The semigroup ring ${\sf k}\left[X_1,\ldots,X_m\right]$ over a field ${\sf k}$
of characteristic 0 associated with ${\sf S}\left({\bf d}^m\right)$ is a
polynomial ring graded by $\deg X_i=d_i,1\leq i\leq m$ and generated by all
monomials $z^{d_i}$. The Hilbert series $H({\bf d}^m;z)$ of a graded ring
${\sf k}\left[z^{d_1},\ldots,z^{d_m}\right]$ is defined by \cite{stan96}
\begin{equation}
H({\bf d}^m;z)=\sum_{s\;\in\;{\sf S}\left({\bf d}^m\right)}z^s=\frac{Q({\bf d}
^m;z)}{\prod_{j=1}^m\left(1-z^{d_j}\right)}\;,\;\;\;z<1\;,\label{hilb0}
\end{equation}
where $Q({\bf d}^m;z)$ is a polynomial in $z$. Thus, we have
\begin{eqnarray}
H({\bf d}^m;z)+\Phi({\bf d}^m;z)=\frac{1}{1-z}\;,\;\;\;z<1\;.\label{nonrep1}
\end{eqnarray}
The number $F\left({\bf d}^m\right):=-1+c\left({\bf d}^m\right)$ is referred 
to as {\em Frobenius number}. As follows from (\ref{def1})
\begin{eqnarray} 
F\left({\bf d}^m\right)=\deg \Phi\left({\bf d}^m;z\right)\;,\;\;\;\;
G\left({\bf d}^m\right)=\sum_{s\;\in\;\Delta\left({\bf d}^m\right)} 1^s=
\Phi\left({\bf d}^m;1\right)\;.\label{frob0}
\end{eqnarray}
The determination of $F\left({\bf d}^m\right)$, $G\left({\bf d}^m\right)$ 
and $H({\bf d}^m;z)$ is called the m--dimensional (mD) Frobenius problem. 
The last restriction ($m\geq 2$) in (\ref{definn0}) is essential since one 
can show that in the one dimensional case the problem is trivial
\begin{eqnarray}
H(d;z)=\frac{1}{1-z^{d}}\;,\;\;\;\;F(d)=G(d)=\infty\;,\;\;\;\;
d\geq 2\;.\label{frob0a}
\end{eqnarray}
The first non--trivial case (m=2) was studied already by J. Sylvester 
\cite{sylv84}     
\begin{equation}
H({\bf d}^2;z)=\frac{1-z^{d_1d_2}}{(1-z^{d_1})(1-z^{d_2})}\;,\;\;\;\;
F({\bf d}^2)=d_1 d_2-d_1-d_2\;,\;\;\;\;G({\bf d}^2)=
\frac{1}{2}(d_1-1)(d_2-1)\;.
\label{sylves1}
\end{equation}
The next non--trivial case (m=3) was extensively studied in the contents of 
commutative 
algebra \cite{herz70}--\cite{denh03} and algebraic geometry of monomial curves 
\cite{kraf85} where a partial progress was achieved (without calculating the 
Hilbert series). A new diagrammatic procedure of construction of the set $\Delta
({\bf d}^3)$ was developed recently in \cite{fel04}. It has paved the way to 
calculate $Q\left({\bf d}^3;z\right)$ and, in accordance with (\ref{nonrep1}) 
and (\ref{frob0}), led to the complete solution of the 3D Frobenius problem. 
Based on Brauer's lemma \cite{brau42} on the matrix representation of the set 
$\Delta({\bf d}^2)$ and Johnson's theorem \cite{john60} on the minimal 
relations, the author \cite{fel04} was able to find the Frobenius number $F
\left({\bf d}^3\right)$, genus $G\left({\bf d}^3\right)$ and the Hilbert 
series $H({\bf d}^3;z)$ of a graded ring for the non--symmetric and symmetric
semigroups ${\sf S}\left({\bf d}^3\right)$ :
\begin{eqnarray}
F\left({\bf d}^3\right)&=&\frac{1}{2}\left[\langle{\bf a},{\bf d}\rangle+
J\left({\bf d}^3\right)\right]-\sum_{i=k}^3d_k,\;\;
G\left({\bf d}^3\right)=\frac{1}{2}\left(1+\langle{\bf a},{\bf d}\rangle-
\sum_{k=1}^3d_k-\prod_{k=1}^3a_{kk}\right)\;,\label{frob1}\\
Q_n({\bf d}^3;z)&=&1-\sum_{k=1}^3 z^{a_{kk}d_k}+z^{1/2\left[\langle{\bf a},
{\bf d}\rangle-J\left({\bf d}^3\right)\right]}+z^{1/2\left[\langle{\bf a},
{\bf d}\rangle+J\left({\bf d}^3\right)\right]}\;,\label{frob2}\\
J^2\left({\bf d}^3\right)&=&\langle{\bf a},{\bf d}\rangle^2-4\sum_{i>j}^3
a_{kk}a_{jj}d_kd_j+4d_1d_2d_3\;,\;\;\langle{\bf a},{\bf d}\rangle=
\sum_{k=1}^3a_{kk}d_k\;,\label{frob3}\\
Q_s({\bf d}^3;z)&=&\left(1-z^{a_{22}d_2}\right)\left(1-z^{a_{33}d_3}\right)\;,
\;\;\;\mbox{if}\;\;\;a_{11}d_1=a_{22}d_2\;,\label{frob4}
\end{eqnarray}
where subscripts "$n$" in (\ref{frob2})  and "$s$" in (\ref{frob4}) stand for 
non--symmetric and symmetric semigroups \cite{fel04}, respectively. The matrix 
$((a_{ij}))$ was introduced by Johnson \cite{john60}. Its elements uniquely 
define the minimal relations for given $d_1,d_2,d_3$
\begin{eqnarray}
&&a_{11}d_1=a_{12}d_2+a_{13}d_3\;,\;\;\;a_{22}d_2=a_{21}d_1+a_{23}d_3\;,\;\;\;
a_{33}d_3=a_{31}d_1+a_{32}d_2\;,\;\;\;\;\mbox{where}\;\;\;\;\;\label{joh1}\\
&&a_{11}=\min\left\{v_{11}\;\bracevert\;v_{11}\geq 2,\;v_{11}d_1=v_{12}d_2+
v_{13}d_3,\;v_{12},v_{13}\in {\mathbb N}\cup\{0\}\right\}\;,\nonumber\\
&&a_{22}=\min\left\{v_{22}\;\bracevert\;v_{22}\geq 2,\;v_{22}d_2=v_{21}d_1+
v_{23}d_3,\;v_{21},v_{23}\in {\mathbb N}\cup\{0\}\right\}\;,\label{joh2}\\
&&a_{33}=\min\left\{v_{33}\;\bracevert\;v_{33}\geq 2,\;v_{33}d_3=v_{31}d_1+
v_{32}d_2,\;v_{31},v_{32}\in {\mathbb N}\cup\{0\}\right\}\;,
\;\;\;\;\mbox{and}\nonumber\\
&&\gcd(a_{11},a_{12},a_{13})=1\;,\;\;\;\gcd(a_{21},a_{22},a_{23})=1\;,\;\;\;
\gcd(a_{31},a_{32},a_{33})=1\;.\nonumber
\end{eqnarray}
Notice that in formulas (\ref{frob1}) -- (\ref{frob4}) only the diagonal 
elements of the Johnson's matrix (\ref{joh2}) appear. Also notice that formulas 
(\ref{frob1}) -- (\ref{frob4}) contain algebraic functions of $a_{jj}$ and 
$d_j,1\leq j\leq 3$. However the $a_{jj}$ cannot be algebraic functions of 
$d_j$ because of the following theorem of Curtis \cite{curt90} :

{\em There is no finite set of polynomials with integer coefficients 
$$
\{f_j(x_1,\ldots,x_m)\}\;,\;\;j=1,\ldots,n
$$

such that for each choice of ${\bf d}^{m},m\geq 3$, there is $j$ such that 
$f_j(d_1,\ldots,d_m)=F({\bf d}^{m})$.}

On the other hand, it would be pretty interesting to build such representations 
for $F\left({\bf d}^3\right)$, $G\left({\bf d}^3\right)$ and $H({\bf d}^3;z)$ 
in the frameworks of analytic number theory avoiding any auxiliary invariants 
like the minimal relations (\ref{joh1}) and (\ref{joh2}). 

Our main result is the analytic representation for the diagonal elements $a_{k
k}\left({\bf d}^3\right)$ of the Johnson's matrix (\ref{joh2}) in terms of 
$d_1$, $d_2$ and $d_3$ :
\begin{eqnarray}
a_{kk}\left({\bf d}^3\right)=\frac{1}{d_k}\lim_{z\to 0}\frac{\ln 
\Psi_k\left({\bf d}^3;z\right)}
{\ln z},\;\;\;\Psi_k\left({\bf d}^3;z\right)=\frac{1}{2\pi}\int_0^{2\pi}
\frac{\left(1-e^{id_jd_lt}\right)dt}{\left(1-e^{id_jt}\right)\left(1-
e^{id_lt}\right)\left(1-z^{d_k}e^{-id_kt}\right)}-1\;,\nonumber \\
\label{joh3}
\end{eqnarray}
where $(k,j,l)=(1,2,3),\;(2,3,1),\;(3,1,2)$.

By (\ref{frob1}) -- (\ref{frob4}) this gives the Frobenius number $F\left({\bf 
d}^3\right)$, genus $G\left({\bf d}^3\right)$ and the Hilbert series $H({\bf d}
^3;z)$, and due to (\ref{nonrep1}) it leads to the generating function 
$\Phi\left({\bf d}^3;z\right)$ for both symmetric and non--symmetric 
semigroups ${\sf S}\left({\bf d}^3\right)$. 
\section{Matrix representation of the set $\Delta\left({\bf d}^2\right)$ and 
the map $\tau$}\label{sect2}
\begin{figure}[h]%[t]
\centerline{\psfig{figure=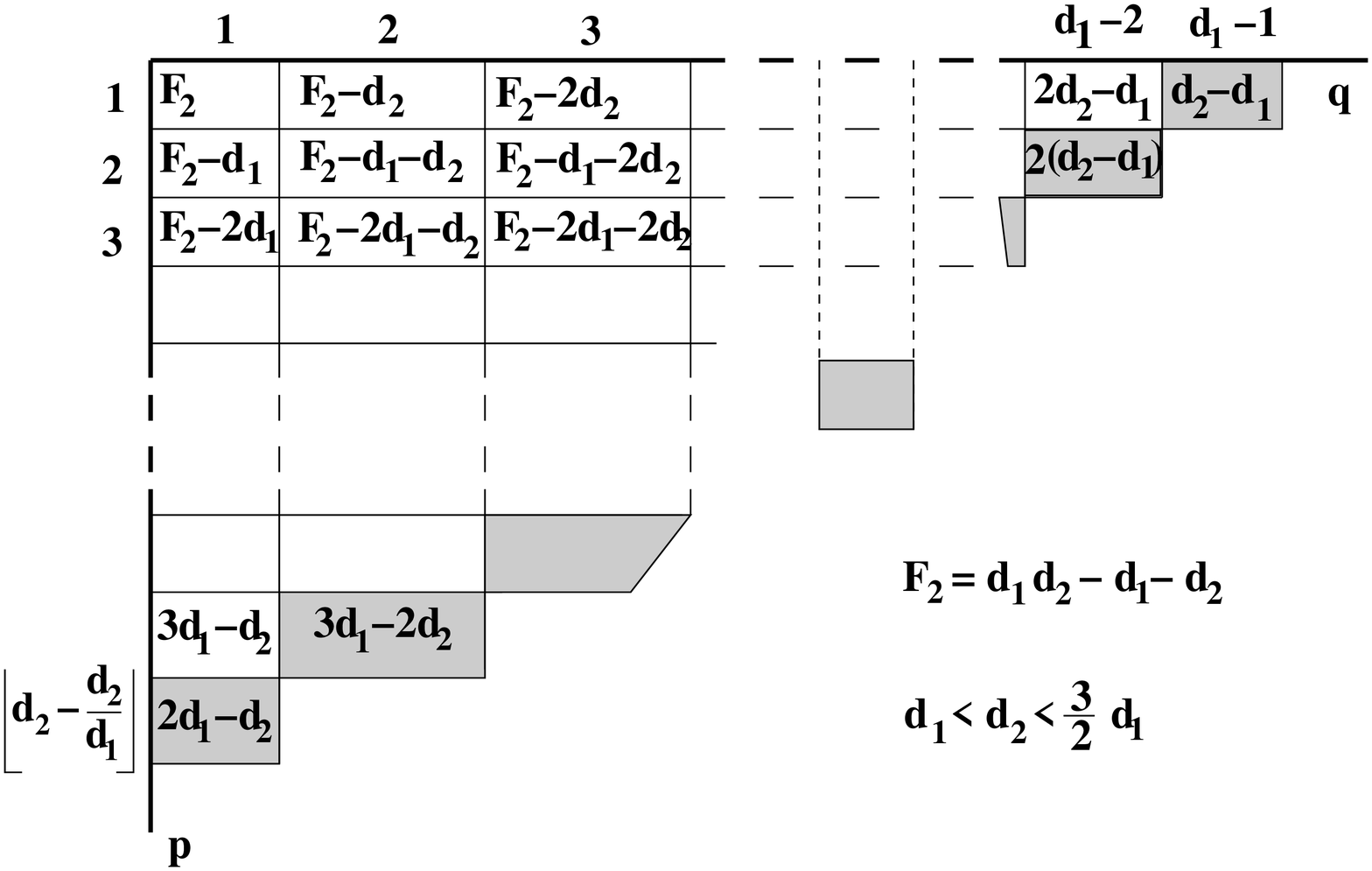,height=5.5cm,width=16cm}}
%width=4.8in}}
\caption{Typical matrix representation $M\left\{\Delta(d_1,d_2)\right\}$ of
the set $\Delta(d_1,d_2)$. The lowest cells in every column of $M\left\{\Delta
(d_1,d_2)\right\}$ ({\em gray color}) is occupied exclusively by the integers 
$1,\ldots,d_1-1$ not in a necessarily consecutive order \cite{fel04}.}
\label{repr1} 
\end{figure}

First, recall the main statements about matrix representation of the set 
$\Delta\left({\bf d}^2\right)$.
\begin{lemma}{\rm (\cite{brau42})}
\label{alem1}
Let $d_1$ and $d_2$ be relatively prime positive integers. Then every positive
integer $s$ not divisible by $d_1$ or by $d_2$ is representable either in the
form $s=xd_1+yd_2,\;x>0,y>0$ or in the form $s=d_1d_2-pd_1-qd_2,\;p,q\in {\mathbb N}$.
\end{lemma}
\begin{definition}
\label{definit1}
Let $2<d_1<d_2$ with $\gcd(d_1,d_2)=1$. Define function $\sigma(p,q)$ as 
follows
\begin{eqnarray}
\sigma(p,q):=d_1 d_2-pd_1-qd_2\;,\;\;\;p,q\in {\mathbb N}\;.\label{syl1}
\end{eqnarray}
\end{definition}
\begin{lemma}{\rm (\cite{fel04})}
\label{alem2}
Let $t$ be an integer and $d_2>d_1$. Then $t\in\Delta\left(d_1,d_2\right)$ iff
$t$ is uniquely representable by
\begin{eqnarray}
t=\sigma(p,q)\;,\;\;\;\;\mbox{where}\label{sylv2}
\end{eqnarray}  
\begin{eqnarray}
1\leq p\leq \left\lfloor d_2-\frac{d_2}{d_1}\right\rfloor\;\;\;\;
1\leq q\leq d_1-1\;,\;\;\;\;\mbox{and}\;\;\;\;
d_1-1\leq \left\lfloor d_2-\frac{d_2}{d_1}\right\rfloor\;,\label{sylv3}
\end{eqnarray}
where $\left\lfloor c\right\rfloor$ denotes an integer part of the number $c$.
\end{lemma}
\begin{lemma}{\rm (\cite{fel04})}
\label{alem3}
Let $2<d_1<d_2$ and $\gcd(d_1,d_2)=1$. Every integer $d_3\in\Delta({\bf d}^2)$
gives rise to the minimal generating set $\{d_1,d_2,d_3\}$ for 
the semigroup ${\sf S}\left({\bf d}^3\right)$ such that $F({\bf d}^3)<
F({\bf d}^2)$.
\end{lemma}
The representation (\ref{sylv2}) of all integers $\sigma(p,q)\in\Delta(d_1,d_2)$
is called {\em the matrix representation} of the set $\Delta(d_1,d_2)$ and is
denoted by $M\left\{\Delta(d_1,d_2)\right\}$ (see Figure \ref{repr1}), so
\begin{eqnarray}
\sigma\left\{M\left\{\Delta(d_1,d_2)\right\}\right\}=\Delta(d_1,d_2)\;,
\label{sylv4}
\end{eqnarray}
i.e. $\sigma$ maps the entries of $M\left\{\Delta(d_1,d_2)\right\}$ onto 
$\Delta(d_1,d_2)$ in a bijective manner (Lemma \ref{alem2}). 

$\sigma(p,q)$ is the integer which occurs in the row $p$ and the column $q$ of
$M\left\{\Delta(d_1,d_2)\right\}$, e.g. $\sigma(1,1)=d_1d_2-d_1-d_2$. Note
that the latter coincides with the Frobenius number $F(d_1,d_2)$.

Throughout the paper we will make use of another entity, {\em a map} $\tau$, 
which was introduced in \cite{fel04} in order to constitute the relationship 
between the set $\Delta\left({\bf d}^3\right)$ and the generating function 
$\Phi\left({\bf d}^3;z\right)$. Actually, the map $\tau$ serves for such 
relationships in any dimension $m\geq 1$ as well.
\begin{definition}
\label{definit2}
The function $\tau$ maps each power series (polynomials including) 
$$
\sum_{s\in {\mathbb N}}c_s\left[\Delta\right]z^s\in{\mathbb N}[z]
$$ 
with $c_s\left[\Delta\right]\in\{0,1\}$ onto the set $\Delta$ of degrees 
$\{s\in \Delta\subset{\mathbb N}\;\bracevert\;c_s\left[\Delta\right]\neq 0\}$.
\end{definition}
The coefficient $c_s\left[\Delta\right]$ is {\em a characteristic function for 
the set} $\Delta\subset {\mathbb N}$
\begin{eqnarray}
c_s\left[\Delta\right]=\left\{\begin{array}{l}1\;,\;\mbox{if}\;\;s\in \Delta
\;,\\0\;,\;\mbox{if}\;\;s\not\in \Delta\;,\end{array}\right.\label{map1bz}   
\end{eqnarray}
which satisfies the following properties :

{\sl Let two sets, $\Delta_1$ and $\Delta_2$, be given such that $\Delta_1,
\Delta_2\subset{\mathbb N}$. Then for set intersections and set unions, we have 
} \cite{Foll99}
\begin{eqnarray}
c_s\left[\Delta_1\cap\Delta_2\right]=c_s\left[\Delta_1\right]c_s\left[\Delta_
1\right]\;,\;\;c_s\left[\Delta_1\cup\Delta_2\right]=c_s\left[\Delta_1\right]+
c_s\left[\Delta_1\right]-c_s\left[\Delta_1\cap\Delta_2\right]\;.\label{map1by}
\end{eqnarray}

\noindent
It also follows from Definition \ref{definit2} and from representations 
(\ref{def1}) and (\ref{hilb0}) that
\begin{eqnarray}
&&\tau\left[\Phi\left({\bf d}^m;z\right)\right]=\Delta\left({\bf d}^m\right)
\;,\;\;\;\tau\left[H\left({\bf d}^m;z\right)\right]={\sf S}\left({\bf d}^m
\right)\;,\;\;\;\tau\left[0\right]=\emptyset\;,\label{map1a}\\
&&\tau^{-1}\left[\Delta\left({\bf d}^m\right)\right]=\Phi\left({\bf d}^m;z
\right)\;,\;\;\;\tau^{-1}\left[{\sf S}\left({\bf d}^m\right)\right]=H\left(   
{\bf d}^m;z\right)\;,\;\;\;\tau^{-1}\left[\emptyset\right]=0\;,\label{map1b} 
\end{eqnarray}  
where $\emptyset$ denotes an empty set. Observe that since all coefficients 
of the polynomial $\Phi\left({\bf d}^m;z\right)$ are 1 or 0, we can uniquely 
reconstruct $\Phi\left({\bf d}^m;z\right)$ from $\Delta\left({\bf d}^m\right)$
and vice versa. In this sense $\tau$ is an isomorphic map. 
The map $\tau$ is also linear in the following sense :
\begin{lemma}
\label{alem4}
Let ${\bf d}^m$ be given and a set $\Delta\left({\bf d}^m\right)$ be related to 
its generating function $\Phi\left({\bf d}^m;z\right)$  by the isomorphic map 
$\tau$ defined in (\ref{map1a}). Let $\Delta_1$ and $\Delta_2$, be subsets 
$\Delta\left({\bf d}^m\right)$. Then the following holds
\begin{eqnarray}
\tau^{-1}\left[\Delta_1\cup \Delta_2\right]=\tau^{-1}\left[\Delta_1\right]+  
\tau^{-1}\left[\Delta_2\right]-\tau^{-1}\left[\Delta_1\cap\Delta_2\right]\;.
\label{map1ba}
\end{eqnarray}
\end{lemma}
{\sf Proof} $\;\;\;$By (\ref{def1}), (\ref{map1bz}), (\ref {map1by}) and 
Definition \ref{definit2} we have
\begin{eqnarray}
\tau^{-1}\left[\Delta_1\cup \Delta_2\right]&=&\sum_{s\;\in\;{\mathbb N}}c_s
\left[\Delta_1\cup \Delta_2\right]z^s=\sum_{s\;\in\;{\mathbb N}}
\left(c_s\left[\Delta_1\right]+c_s\left[\Delta_1\right]-c_s\left[\Delta_1
\cap\Delta_2\right]\right)z^s=\nonumber\\
&&\tau^{-1}\left[\Delta_1\right]+\tau^{-1}\left[
\Delta_2\right]-\tau^{-1}\left[\Delta_1\cap\Delta_2\right]\;,\nonumber
\end{eqnarray}
that proves the Lemma.$\;\;\;\;\;\;\Box$

The following statement is a consequence of Lemma \ref{alem4} and generalizes 
the result obtained in \cite{fel04} (formula (64))
\begin{lemma} 
\label{alem5}
Let ${\bf d}^m$ be given and a set $\Delta\left({\bf d}^m\right)$ be related 
to its generating function $\Phi\left({\bf d}^m;z\right)$  by the map $\tau$ 
according to (\ref{map1a}) and (\ref{map1b}). Let $n$ sets $\Delta_i,\;i=1,
\ldots,n$ be given such that
\begin{eqnarray}
\Delta_i\subset\Delta\left({\bf d}^m\right)\;,\;\;\;\Delta_i\cap\Delta_j=
\emptyset\;,\;\;\;\;\mbox{for}\;\;\;\;i\neq j=1,\ldots,n\;.\label{amap1} 
\end{eqnarray} 
Then the following holds
\begin{eqnarray}
\tau^{-1}\left[\bigcup_{i=1}^n\Delta_i\right]=\sum_{i=1}^n\tau^{-1}\left[
\Delta_i\right]\;.\label{amap2}
\end{eqnarray}
\end{lemma}
{\sf Proof} $\;\;\;$By induction on $n$ applying Lemma \ref{alem4} to the 
left hand side of (\ref{amap2}) consecutively and making use of (\ref{amap1}) 
we obtain
\begin{eqnarray}
\tau^{-1}\left[\bigcup_{i=1}^n \Delta_i\right]=\tau^{-1}\left[\bigcup_{i=1}
^{n-1}\Delta_i\right]+\tau^{-1}\left[\Delta_n\right]-\tau^{-1}\left[\left\{
\bigcup_{i=1}^{n-1}\Delta_i\right\}\bigcap\Delta_n\right]=\nonumber\\
\tau^{-1}\left[\bigcup_{i=1}^{n-1}\Delta_i\right]+\tau^{-1}\left[\Delta_n
\right]-\tau^{-1}\left[\bigcup_{i=1}^{n-1}\left\{\Delta_i\cap\Delta_n\right\}
\right]=\tau^{-1}\left[\bigcup_{i=1}^{n-1}\Delta_i\right]+\tau^{-1}
\left[\Delta_n\right]=\ldots=\sum_{i=1}^n\tau^{-1}\left[\Delta_i\right]\nonumber
\end{eqnarray}
that proves the Lemma.$\;\;\;\;\;\;\Box$

\subsection{The intersection set $\Delta\left(d_1,d_2\right)\cap{\sf 
S}\left(d_3\right)$ and its generating function}\label{sect21}
In this Section we construct the intersection set $\Delta\left(d_1,d_2\right)
\cap{\sf S}\left(d_3\right)$. This set is essential to determine an 
explicit non--algebraic expression for $a_{33}(d_1,d_2,d_3)$ (see next Section).

First, define a 1D numerical semigroup ${\sf S}\left(d_3\right)$ of 
integers $\sigma=0 \bmod (d_3)$
\begin{eqnarray}
{\sf S}\left(d_3\right):=\left\{\sigma\;\bracevert\;\sigma=jd_3,\;j\in{\mathbb 
N}\cup\{0\},\;d_3\geq 2\right\}\;,\label{inter1}
\end{eqnarray}
which is generated by the Hilbert series $H\left(d_3;z\right)$ given by 
(\ref{frob0a}). 

Consider the intersection set $\Delta\left(d_1,d_2\right)\cap{\sf S}\left(d_3
\right)$ which consists exclusively of the integers $\sigma=jd_3$, the index 
$j$ runs with jumps from 1 to $N_3$ such that $N_3d_3\in\Delta\left(d_1,d_2
\right)$ still holds. Notice that $N_3$ satisfies
$$
N_3\leq \frac{d_1 d_2-d_1-d_2}{d_3}< d_1-1\;,
$$
that follows from $N_3d_3\leq F(d_1,d_2)$ and $d_2(d_1-1)-d_1<d_3(d_1-1)$.

\begin{figure}[h]%[t]
\centerline{\psfig{figure=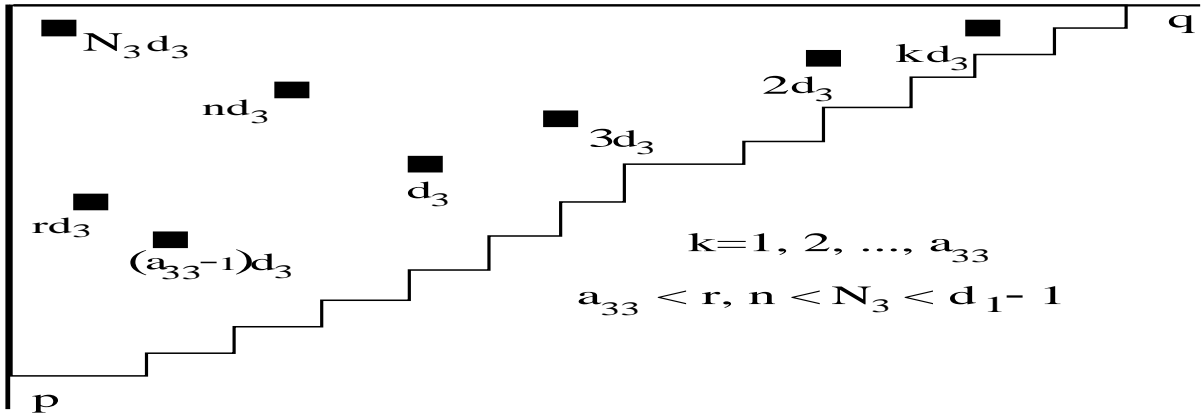,height=5.5cm,width=16cm}}
\caption{Integers $jd_3\in\Delta\left(d_1,d_2\right)\cap{\sf S}\left(d_3\right)$
({\em black boxes}) are distributed inside the matrix representation $M\left\{
\Delta(d_1,d_2)\right\}$.}
\label{repr22}
\end{figure}

The integers $jd_3$ are distributed inside the matrix representation $M\left\{
\Delta(d_1,d_2)\right\}$ as shown in Figure \ref{repr22}. Let us prove an 
important Lemma. 

\begin{lemma}
\label{alem6}
Let ${\bf d}^3$ be given, ${\bf d}^3=(d_1,d_2,d_3)$, and assume that the 
minimal relations are defined by (\ref{joh1}) and (\ref{joh2}). Then
\begin{eqnarray}
jd_3\in\Delta\left(d_1,d_2\right)\cap{\sf S}\left(d_3\right)\;,\;\;\;j=1,
\ldots,a_{33}-1\;,\label{inter3}
\end{eqnarray}
and 
\begin{eqnarray}
a_{33}d_3\not\in\Delta\left(d_1,d_2\right)\cap{\sf S}\left(d_3\right)\;.
\label{inter3a}
\end{eqnarray}
\end{lemma}
Notice that $j$ runs in (\ref{inter3}) without jumps.

\noindent
{\sf Proof} $\;\;\;$The proof follows from two obvious relations. First, 
by (\ref{inter1}) we have
\begin{eqnarray}
jd_3\in {\sf S}\left(d_3\right)\;,\;\;\;j=1,\ldots,a_{33},\ldots\;,
\label{inter4}
\end{eqnarray}
and, next, by the 3rd line of (\ref{joh2}) we have
\begin{eqnarray}
jd_3\in\Delta\left(d_1,d_2\right)\;,\;\;j=1,\ldots,a_{33}-1\;,\;\;\mbox{and}
\;\;\;a_{33}d_3\not\in\Delta\left(d_1,d_2\right)\;,
\label{inter5}
\end{eqnarray}
that leads to (\ref{inter3}) and (\ref{inter3a}).$\;\;\;\;\;\;\Box$

\noindent
Keeping in mind Lemma \ref{alem6} we write the generating function of the set 
$\Delta\left(d_1,d_2\right)\cap{\sf S}\left(d_3\right)$. 

According to (\ref{map1b}) we get
\begin{eqnarray}
\tau^{-1}\left[\Delta\left(d_1,d_2\right)\cap{\sf S}\left(d_3\right)\right]=
\sum_{j=1}^{a_{33}-1}z^{jd_3}+\sum_{j>a_{33}}^{N_3}z^{jd_3}\;,\label{inter6}
\end{eqnarray}
where index $j$ runs with jumps in the second summation of (\ref{inter6}). 
Define an auxiliary function $\Psi_3\left({\bf d}^3;z\right)$ by
\begin{eqnarray}
\Psi_3\left({\bf d}^3;z\right)=\frac{z^{d_3}}{1-z^{d_3}}-\tau^{-1}\left[\Delta
\left(d_1,d_2\right)\cap{\sf S}\left(d_3\right)\right]\;.\label{inter7}  
\end{eqnarray}
It plays a key role in calculating the diagonal element $a_{33}$. 
\begin{theorem}
\label{theo1}
Let ${\bf d}^3$ be given, ${\bf d}^3=(d_1,d_2,d_3)$, and assume that the 
minimal relations are defined by (\ref{joh1}) and (\ref{joh2}). Then the 
diagonal element $a_{33}\left({\bf d}^3\right)$ is given by
\begin{eqnarray}
a_{33}\left({\bf d}^3\right)=\frac{1}{d_3}\lim_{z\to 0}\frac{\ln 
\Psi_3\left({\bf d}^3;z\right)}{\ln z}\;.\label{inter8}
\end{eqnarray}
\end{theorem}
{\sf Proof} $\;\;\;$A series expansion for $\Psi_3\left({\bf d}^3;z\right)$ 
reads 
\begin{eqnarray}
\Psi_3\left({\bf d}^3;z\right)=z^{a_{33}d_3}+\sum_{j>a_{33}}^{j<N_3}z^{jd_3}+
\sum_{j=N_3+1}^{\infty}z^{jd_3}\;,\label{inter9}
\end{eqnarray}
where index $j$ runs with jumps in the second summation of (\ref{inter9}). The 
power series (\ref{inter9}) converges for $z<1$
$$
\Psi_3\left({\bf d}^3;z\right)=z^{a_{33}d_3}+\sum_{j>a_{33}}^{j<N_3}z^{jd_3}+
\frac{z^{(N_3+1)d_3}}{1-z^{d_3}}\;.
$$
The first term in (\ref{inter9}) turns to be the leading term if $z\to 0$, 
i.e. $\Psi_3\left({\bf d}^3;z\right)\stackrel{z\to 0}\longrightarrow z^{a_{33}
d_3}$. Hence, formula (\ref{inter8}) follows.$\;\;\;\;\;\;\Box$

Formulas (\ref{inter7}) and (\ref{inter8}) have nothing specific concerned with
arrangements of the $d_i$ in the tuple $(d_1,d_2,d_3)$. Therefore the similar 
formulas can be obtained for the diagonal elements $a_{11}\left({\bf d}^3\right)$
and $a_{22}\left({\bf d}^3\right)$
\begin{eqnarray}
a_{kk}\left({\bf d}^3\right)=\frac{1}{d_k}\lim_{z\to 0}\frac{\ln 
\Psi_k\left({\bf d}^3;z\right)}{\ln z}\;,\;\;\;k=1,2\;,\label{inter10}
\end{eqnarray}
where
\begin{eqnarray}
\Psi_1\left({\bf d}^3;z\right)=\frac{z^{d_1}}{1-z^{d_1}}-\tau^{-1}\left[\Delta
\left(d_2,d_3\right)\cap{\sf S}\left(d_1\right)\right]\;,\;\;\;\Psi_2\left({\bf 
d}^3;z\right)=\frac{z^{d_2}}{1-z^{d_2}}-\tau^{-1}\left[\Delta
\left(d_3,d_1\right)\cap{\sf S}\left(d_2\right)\right]\;.\nonumber
\end{eqnarray}  
\section{Analytic representation of generating functions}\label{sect3}
In the previous Section we have found the diagonal elements $a_{jj}\left({\bf 
d}^3\right),j=1,2,3$, through the auxiliary functions $\Psi_j\left({\bf d}^3;
z\right),j=1,2,3$. Thus, the non--algebraic expressions for 
$a_{jj}\left({\bf d}^3\right)$ will be given in the closed 
form if the analytic representation of corresponding generating functions 
$$
\tau^{-1}\left[\Delta\left(d_j,d_l\right)\cap{\sf S}\left(d_k\right)\right]\;,
\;\;\;(k,j,l)=(1,2,3),\;(2,3,1),\;(3,1,2)\;,
$$
will be found. The present Section is completely devoted to this question. In 
order to resolve it we make use of the {\em Hadamard product} of analytic 
functions.
\subsection{Hadamard product of analytic functions}\label{sect31}
Let $u(z)=\sum_{n=0}^{\infty} a_nz^n$ and $v(z)=\sum_{n=0}^{\infty}b_nz^n$ be 
analytic functions in the unit disk $|z|<1$. Define the $\otimes$--product 
of $u(z)$ and $v(z)$ by
\begin{eqnarray}
u(z)\otimes v(z)=\frac{1}{2\pi}\int_0^{2\pi}u\left(ze^{it}\right)v\left(
e^{-it}\right)dt=\sum_{n=0}^{\infty}a_nb_nz^n\;,\;\;\;|z|<1\;.\label{had1}
\end{eqnarray}
This product was introduced by J. Hadamard \cite{hada99} to discuss the 
singularities of the analytic function $f(z)$ defined by the series $\sum_{n=0}
^{\infty}a_nb_nz^n$ in terms of those of the functions $u$ and $v$ ({\em 
"multiplication of singularities"} theorem \cite{titc39}). 

Another representation, which is similar to (\ref{had1}), can be found for the 
analytic function $f(z^2)=\sum_{n=0}^{\infty}a_nb_nz^{2n}$ through the 
integral convolution "$\circ$" of the functions $u$ and $v$,
\begin{eqnarray}
u(z)\circ v(z)=\frac{1}{2\pi}\int_0^{2\pi}u\left(ze^{it}\right)v\left(ze^{-it}
\right)dt=\sum_{n=0}^{\infty}a_nb_nz^{2n}\;,\;\;\;|z|<1\;.\label{had2}
\end{eqnarray}
This product was also used by J. Hadamard \cite{hada99} and is strongly 
related to the $\otimes$--product
\begin{eqnarray}
u\left(z^2\right)\otimes v\left(z^2\right)=u(z)\circ v(z)\;,\label{had3}
\end{eqnarray}
e.g. for the analytic function $u(z)$ we have
\begin{eqnarray}
u(z)\otimes (1-z)^{-1}=u(z)\;,\;\;\;\;\mbox{but}\;\;\;\;\;
u(z)\circ (1-z)^{-1}=u\left(z^2\right)\;.\label{had2a}
\end{eqnarray}
Both $\otimes$-- and $\circ$--products have been variously referred to as {\em 
the Hadamard product}, {\em quasi inner product} or {\em the termwise product} 
(see survey \cite{horn90}). The $\otimes$-- and $\circ$--products appear 
naturally in functional analysis, e.g. the Bieberbach conjecture for univalent 
analytic functions \cite{loew59} and the Polya-Schoenberg conjecture for 
analytic convex mappings \cite{rush73}.

There is a variety of algebraic and other analytic properties associated with 
$\circ$-- and $\otimes$--products \cite{brag86}. The commutative and 
distributive laws can be easily verified
\begin{eqnarray}
&&u(z)\circ v(z)=v(z)\circ u(z)\;,\;\;\;u(z)\circ [v(z)+w(z)]=u(z)\circ v(z)+
u(z)\circ w(z)\;,\nonumber\\
&&u(z)\otimes v(z)=v(z)\otimes u(z)\;,\;\;\;u(z)\otimes [v(z)+w(z)]=u(z)
\otimes v(z)+u(z)\otimes w(z)\;.\label{had3a}
\end{eqnarray}
Notice that (\ref{had2}) manifests the non--associativity of the 
$\circ$--product 
$$
[u(z)\circ v(z)]\circ w(z)\neq u(z)\circ [v(z)]\circ w(z)]\;,
$$ 
e.g. $[z\circ z]\circ z^2=z^4$, but $z\circ [z\circ z^2]=0$. However the 
$\otimes$--product is associative, 
\begin{eqnarray}
[u(z)\otimes v(z)]\otimes w(z)=u(z)\otimes [v(z)]\otimes w(z)]\;,\label{had3b}
\end{eqnarray}
e.g. $[z\otimes z]\otimes z^2=z\otimes [z\otimes z^2]=0$. The last property is 
very important to generalize the Hadamard product. We give an appropriate 
extension of $\otimes$--product referred to as {\em the Hadamard multiple 
product} \cite{brag99}. 

Let $u_k(z)=\sum_{n=0}^{\infty} a_n^{(k)}z^n,k=1,\ldots,N$, be analytic 
functions in the unit disk $|z|<1$. Define the Hadamard multiple product 
$\bigotimes$ of $u_k(z)$ by
\begin{eqnarray}
\bigotimes_{k=1}^Nu_k(z)=u_1(z)\otimes u_2(z)\otimes \ldots\otimes u_N(z)=\sum_
{n=0}^{\infty}a_n^{(1)}\cdot a_n^{(2)}\cdot\ldots\cdot a_n^{(N)}\cdot z^n\;.
\label{had3c}
\end{eqnarray}
The $\bigotimes$--product is well defined due to the associativity of the 
$\otimes$--product (\ref{had3b}). Its integral representation has the form
\begin{eqnarray}
\bigotimes_{k=1}^Nu_k(z)=\frac{1}{(2\pi)^{N-1}}\int_0^{2\pi}\ldots\int_0^{2\pi}
u_1\left(ze^{i\sum_{k=2}^Nt_k}\right)\prod_{k=2}^Nu_k\left(e^{-it_k}\right)dt_k
\;.\label{had3d}
\end{eqnarray}
We give some examples of the $\otimes$--product operating with rational 
functions \cite{titc39}
\begin{eqnarray}
\bigotimes_{k=1}^N\frac{1}{a_k-z}=\frac{1}{\prod_{k=1}^Na_k-z}\;,\;\;\;\;\;\;
\bigotimes_{k=1}^N\frac{1}{1-z^{n_k}}=\frac{1}{1-z^{L}}\;,\;\;\;n_k\in {\mathbb 
N}\;,\;\;\;L={\rm lcm}\left(n_1,\ldots,n_N\right)\label{had4}
\end{eqnarray}
and also of the $\circ$--product operating with transcendental functions  
\cite{brag86}, \cite{brag99}
\begin{eqnarray}
\cosh\frac{z}{2}\circ\cosh\frac{z}{2}+\sinh\frac{z}{2}\circ\sinh\frac{z}{2}=   
I_0(z)\;,\;\;\;\cosh\frac{z}{2}\circ\cosh\frac{z}{2}-\sinh\frac{z}{2}\circ\sinh
\frac{z}{2}=J_0(z)\;,\label{had5}\\
e^{z/2}\circ e^{z/2}=I_0(z)\;,\;\;\;e^{-z}\circ (1+z)^n=L_n(z^2)\;,\;\;\;
\frac{1}{(1-z)^{\alpha}}\circ \frac{1}{(1-z)^{\beta}}=\;_2F_1\left(\alpha;
\beta;1;z^2\right)\;,\label{hada5}
\end{eqnarray}
where $\alpha,\beta>0$. In (\ref{had5}) and (\ref{hada5}) the functions 
$I_0(z)$, $J_0(z)$, $L_n(z)$ and $_2F_1\left(\alpha;\beta;1;z\right)$ denote 
the modified and non--modified Bessel functions of the 1st kind, the Laguerre 
polynomial and the hypergeometric function, respectively. 

Because of its connection with the integral convolution (\ref{had2}), the 
Hadamard product leads to elegant evaluations of complicated trigonometric 
integrals and provide analytic derivations of combinatorial identities 
\cite{brag86}, \cite{brag99}. Finally, the Hadamard product was used in 
\cite{brag94} to construct solutions for a variety of Cauchy--type problems.

It turns out that the $\otimes$--product is a sensitive tool to deal with 
an intersection of discrete numerical sets. This will be the subject for 
discussion in the next Section \ref{sect32}. In order to be technically equiped 
we derive here two important formulas which will be useful further.
\begin{lemma}
\label{alem7}
Let $u(z)=\sum_{k=0}^{\infty} a_kz^k$ be an analytic function in the unit 
disk $|z|<|c|$. Then
\begin{eqnarray}
u(z)\otimes\frac{1}{c-z}=\frac{1}{c}\;u\left(\frac{z}{c}\right)\;,\;\;\;
c\neq 0\;.\label{had5o}
\end{eqnarray}
\end{lemma}
{\sf Proof} $\;\;\;$Calculating the $\otimes$--product we get 
\begin{eqnarray}
u(z)\otimes\frac{1}{c-z}=\frac{1}{c}\;u(z)\otimes\frac{1}{1-z/c}=\frac{1}{c}\;
\sum_{k=0}^{\infty}a_kz^k\otimes\sum_{k=0}^{\infty}\left(\frac{z}{c}\right)^k=
\frac{1}{c}\;\sum_{k=0}^{\infty}a_k\frac{z^k}{c^k}=
\frac{1}{c}\;u\left(\frac{z}{c}\right)\;,\nonumber
\end{eqnarray}
that proves formula (\ref{had5o}).$\;\;\;\;\;\;\Box$

Before going to the next formula (Lemma \ref{alem81}) we prove an 
elementary 
identity
\begin{eqnarray}
\frac{n}{z^n-1}=\sum_{l=1}^n\frac{w_n^l}{z-w_n^l}\;,\;\;\;w_n=\exp\left(
\frac{2\pi i}{n}\right)\;.\label{iden1}
\end{eqnarray}
Rewrite the right hand side of (\ref{iden1}) in the form
\begin{eqnarray}
\sum_{l=1}^n\frac{w_n^l}{z-w_n^l}=\frac{1}{z^n-1}\sum_{k=1}^n(-1)^{k-1}k\Pi_k
z^{n-k}=\frac{\Pi_1z^{n-1}-2\Pi_2z^{n-2}+\ldots+n(-1)^{n-1}\Pi_n}{z^n-1}\;,
\label{iden2}
\end{eqnarray}
where $\Pi_k$ denote the basic symmetric polynomials
\begin{eqnarray}  
\Pi_k=\sum_{l_1>l_2>...>l_k=1}^{n}w_n^{l_1}w_n^{l_2}\ldots w_n^{l_k},\;\;\;
\mbox{i.e.}\;\;\;\Pi_1=\sum_{l=1}^{n}w_n^l,\;\;\Pi_2=\sum_{l_1>l_2=1}^{n}
w_n^{l_1}w_n^{l_2}\;,\;\ldots\;,\;\Pi_n=\prod_{l=1}^{n}w_n^l.\label{iden3}
\end{eqnarray}
Recall Vieta's formula for the sum $S_k$ of the products of distinct roots 
$z_i$ of the polynomial equation of degree $n$  
\begin{eqnarray}
a_nz^n+a_{n-1}z^{n-1}+\ldots+a_1z+a_0=0\;,\nonumber
\end{eqnarray}
which reads \cite{waerd93}
\begin{eqnarray}
S_k=(-1)^k\frac{a_{n-k}}{a_n}\;,\;\;\mbox{where}\;\;S_k=\sum_{l_1>l_2>...>
l_k=1}^{n}z_{l_1}z_{l_2}\ldots z_{l_k}.\label{iden4}
\end{eqnarray}
Since the $w_n^l$--roots are associated with the polynomial equation $z^n-1=0$ 
then the Vieta's formula (\ref{iden4}) gives for $\Pi_k$--polynomials
\begin{eqnarray}
\Pi_k=\left\{\begin{array}{r}0\;,\;\mbox{if}\;\;1\leq k<n\;,\\
(-1)^{n+1}\;,\;\mbox{if}\;\;k=n\;.\end{array}\right.\label{iden6}
\end{eqnarray}
Substituting (\ref{iden6}) into (\ref{iden2}) we arrive at (\ref{iden1}).
\begin{lemma}
\label{alem81}
Let $u(z)=\sum_{k=0}^{\infty} a_kz^k$ be an analytic function in the unit
disk $|z|<1$. Then
\begin{eqnarray}
u(z)\otimes\frac{1}{1-z^n}=\frac{1}{n}\;\sum_{k=0}^{n-1}u\left(zw_n^k\right)\;.
\label{had5p}
\end{eqnarray}
\end{lemma}
{\sf Proof} $\;\;\;$Making use of the identity (\ref{iden1}) we obtain by 
Lemma \ref{alem7}
\begin{eqnarray}
u(z)\otimes\frac{1}{1-z^n}=\frac{1}{n}\;\sum_{k=1}^n\left(u(z)\otimes\frac{
w_n^l}{w_n^l-z}\right)=\frac{1}{n}\sum_{l=1}^n\frac{w_n^l}{w_n^l}u\left(
\frac{z}{w_n^l}\right)=\frac{1}{n}\sum_{l=1}^nu\left(zw_n^{n-l}\right)=
\frac{1}{n}\sum_{k=0}^{n-1}u\left(zw_n^k\right)\nonumber
\end{eqnarray}
that proves formula (\ref{had5p}).$\;\;\;\;\;\;\Box$

Henceforth we use the $\otimes$--product and refer to as the Hadamard product.
\subsection{Representation of $\tau^{-1}\left[\Delta\left(d_1,d_2\right)\cap{
\sf S}\left(d_3\right)\right]$}\label{sect32}
In this Section we derive the analytic representation of the generating 
function for the intersection of two sets, $\Delta\left(d_1,d_2\right)$ and 
${\sf S}\left(d_3\right)$. Start with the following proposition which utilizes 
the Hadamard product for the intersecting sets $\Delta_1$ and $\Delta_2$. 
\begin{lemma} 
\label{alem91}
Let $\Delta_1,\Delta_2$ be subsets of ${\mathbb N}$, and let their 
corresponding generating functions $\tau^{-1}\left[\Delta_1\right]$ and 
$\tau^{-1}\left[\Delta_2\right]$ be
\begin{eqnarray}
\tau^{-1}\left[\Delta_i\right]=\sum_{s\in {\mathbb N}}c_s
\left[\Delta_i\right]z^s\;,\;\;\;i=1,2\;,\label{had5a}
\end{eqnarray}
where $c_s\left[\Delta_i\right]$ stands for corresponding characteristic 
function of the set $\Delta_i$ and satisfies (\ref{map1bz}). Then the set 
$\Delta_1\cap\Delta_2$ is generated by
\begin{eqnarray}
\tau^{-1}\left[\Delta_1\cap\Delta_2\right]=\tau^{-1}\left[\Delta_1\right]
\otimes\tau^{-1}\left[\Delta_2\right]\;.\label{had5b}
\end{eqnarray}
\end{lemma}
{\sf Proof} $\;\;\;$By Definition \ref{definit2} of generating function 
$\tau^{-1}\left[\Delta\right]$ and the 1st relation in (\ref{map1by}) we 
have
\begin{eqnarray}
\tau^{-1}\left[\Delta_1\cap\Delta_2\right]=\sum_{s\in {\mathbb N}}c_s\left[
\Delta_1\cap\Delta_2\right]z^s=\sum_{s\in {\mathbb N}}c_s\left[\Delta_1\right]
c_s\left[\Delta_1\right]z^s\;.\label{had5c}
\end{eqnarray}
On the other hand, by definition (\ref{had1}) of the Hadamard product we have
\begin{eqnarray}
\tau^{-1}\left[\Delta_1\right]\otimes\tau^{-1}\left[\Delta_2\right]=
\sum_{s\in {\mathbb N}}c_s\left[\Delta_1\right]z^s\otimes
\sum_{s\in {\mathbb N}}c_s\left[\Delta_2\right]z^s=\sum_{s\in {\mathbb N}}c_s
\left[\Delta_1\right]c_s\left[\Delta_1\right]z^s\;.\label{had5d}
\end{eqnarray}
A comparison of (\ref{had5c}) and (\ref{had5d}) proves the Lemma.
$\;\;\;\;\;\;\Box$

Lemma \ref{alem91} has a simple generalization
\begin{eqnarray}
\tau^{-1}\left[\cap_{j=1}^N\Delta_j\right]=\bigotimes_{j=1}^N\tau^{-1}\left[
\Delta_j\right]\;.\label{had5dd}
\end{eqnarray}
Returning to $\tau^{-1}\left[\Delta\left(d_1,d_2\right)\cap{\sf S}\left(d_3
\right)\right]$ we can verify that both generating functions $\tau^{-1}\left[
\Delta\left(d_1,d_2\right)\right]$ and $\tau^{-1}\left[{\sf S}\left(d_3
\right)\right]$ are representable in the form (\ref{had5a})
\begin{eqnarray}
\tau^{-1}\left[\Delta\left(d_1,d_2\right)\right]=\Phi\left(\{d_1,d_2\};z\right)
=\sum_{s\;\in\;\Delta\left(d_1,d_2\right)}z^s\;,\;\;\;\;\;\;\tau^{-1}\left[{\sf 
S}\left(d_3\right)\right]=\frac{1}{1-z^{d_3}}=\sum_{s\;\in\;{\sf S}\left(d_3
\right)}z^s\;.\nonumber
\end{eqnarray}
Therefore we obtain
\begin{eqnarray}
\tau^{-1}\left[\Delta\left(d_1,d_2\right)\cap{\sf S}\left(d_3\right)\right]=
\tau^{-1}\left[\Delta\left(d_1,d_2\right)\right]\otimes\tau^{-1}\left[{\sf 
S}\left(d_3\right)\right]=\Phi\left(\{d_1,d_2\};z\right)\otimes 
\frac{1}{1-z^{d_3}}\;.\label{had5f}
\end{eqnarray}
Formula (\ref{had5f}) can be slightly simplified by utilizing the relationship 
(\ref{nonrep1}) between the generating function $\Phi$ and the Hilbert series 
$H$ and of distributive law (\ref{had3a}) for the Hadamard product
\begin{eqnarray}
\tau^{-1}\left[\Delta\left(d_1,d_2\right)\cap{\sf S}\left(d_3\right)\right]=
\frac{1}{1-z}\otimes\frac{1}{1-z^{d_3}}-H_{12}(z)\otimes 
\frac{1}{1-z^{d_3}}\;,\label{had6}
\end{eqnarray}
where $H_{12}(z)$ denotes for short the Hilbert series $H\left(\{d_1,d_2\};z
\right)$. On the last step we make use of the 1st equality in (\ref{had2a}) and 
get the generating function for the intersection of two sets, $\Delta\left(d_1,
d_2\right)$ and ${\sf S}\left(d_3\right)$,
\begin{eqnarray}
\tau^{-1}\left[\Delta\left(d_1,d_2\right)\cap{\sf S}\left(d_3\right)\right]=
\frac{1}{1-z^{d_3}}-H_{12}(z)\otimes\frac{1}{1-z^{d_3}}\;,
\label{had6a}
\end{eqnarray}
and also the function $\Psi_3\left({\bf d}^3;z\right)$ introduced in 
(\ref{inter7})
\begin{eqnarray}
\Psi_3\left({\bf d}^3;z\right)=\frac{z^{d_3}}{1-z^{d_3}}-\tau^{-1}\left[\Delta
\left(d_1,d_2\right)\cap{\sf S}\left(d_3\right)\right]=H_{12}(z)\otimes\frac{1}
{1-z^{d_3}}-1\;.\label{had7}
\end{eqnarray}
Applying now Lemma \ref{alem81} to formula (\ref{had7}) we get
\begin{eqnarray}
\Psi_3\left({\bf d}^3;z\right)=\frac{1}{d_3}\;\sum_{k=0}^{d_3-1}H_{12}\left(z
w_{d_3}^k\right)-1\;,\;\;\;w_{d_3}=\exp\left(\frac{2\pi i}{d_3}\right)\;.
\label{had7x}
\end{eqnarray}
Simplifying the expression (\ref{sylves1}) for $H_{12}\left(z\right)$
\begin{eqnarray}
H_{12}\left(z\right)=\sum_{p=0}^{d_2-1}\frac{z^{pd_1}}{1-z^{d_2}}=
\sum_{p=0}^{d_2-1}\sum_{q=0}^{\infty}z^{pd_1+qd_2}\;,\label{had7z}
\end{eqnarray}
and substituting it into (\ref{had7x}) we obtain
\begin{eqnarray}
\Psi_3\left({\bf d}^3;z\right)=\frac{1}{d_3}\;\sum_{p=0}^{d_2-1}\sum_{q=0}^
{\infty}z^{pd_1+qd_2}\sum_{k=0}^{d_3-1}w_{d_3}^{k(pd_1+qd_2)}-1\;.\label{had7y}
\end{eqnarray}   
The inner sum in (\ref{had7y}) does vanish for $p,q$ such that $pd_1+qd_2$ 
is not divisible by $d_3$. Indeed, 
\begin{eqnarray}
\sum_{k=0}^{d_3-1}w_{d_3}^{k(pd_1+qd_2)}=\frac{\exp\left[2\pi i\;(pd_1+qd_2)
\right]-1}{\exp\left(2\pi i\;\frac{pd_1+qd_2}{d_3}\right)-1}=0\;,\;\;\;
\mbox{if}\;\;\;d_3\not|\;pd_1+qd_2\;.\nonumber
\end{eqnarray}
Thus, it remains
\begin{eqnarray} 
\Psi_3\left({\bf d}^3;z\right)=\sum_{p=0}^{d_2-1}\sum_{q=0\atop d_3\mid\;pd_1+
qd_2}^{\infty}z^{pd_1+qd_2}-1=\sum_{p=1}^{d_2-1}\sum_{q=1\atop d_3\mid\;pd_1+
qd_2}^{\infty}z^{pd_1+qd_2}=\sum_{j=1}^{\infty}z^{b_jd_3}\;,\label{had7w}
\end{eqnarray}
where $b_j\in {\mathbb N}$ is defined as the integer which satisfies the 
Diophantine equations in $b_j$
\begin{eqnarray}
pd_1+qd_2=b_jd_3\;,\;\;\;p=0,\ldots,d_2-1\;,\;\;\;q=0,\ldots\label{had7v}
\end{eqnarray}
at least with one solution. Recalling the definition (\ref{joh1}) and 
(\ref{joh2}) of the Johnson's minimal relations we conclude that 
\begin{eqnarray}
b_1=a_{33}\;,\label{had777}
\end{eqnarray}
and the first term in series expansion (\ref{had7w}) reads $z^{a_{33}d_3}$.   
Hence, formula (\ref{inter8}) follows. Notice that (\ref{had7w}) can 
be obtained by straightforward calculation of the Hadamard product (\ref{had7})
according to Lemma \ref{alem91} and representation (\ref{had7z})
\begin{eqnarray}
\Psi_3\left({\bf d}^3;z\right)=\left(\sum_{p=0}^{d_2-1}\sum_{q=0}^{
\infty}z^{pd_1+qd_2}\right)\otimes \sum_{r=0}^{\infty}z^{rd_3}-1=
\left(\sum_{p=1}^{d_2-1}\sum_{q=1}^{\infty}z^{pd_1+qd_2}\right)\otimes 
\sum_{r=1}^{\infty}z^{rd_3}=\sum_{j=1}^{\infty}z^{b_jd_3}\;,\nonumber
\end{eqnarray}
where $b_j$ is defined in (\ref{had7v}). 

We finish this Section by giving an integral representation for (\ref{had7}) 
according to (\ref{had1})
\begin{eqnarray}
\Psi_3\left({\bf d}^3;z\right)=\frac{1}{2\pi}\int_0^{2\pi}\frac{\left(1-
e^{id_1d_2t}\right)dt}{\left(1-e^{id_1t}\right)\left(1-e^{id_2t}\right)\left(
1-z^{d_3}e^{-id_3t}\right)}-1\;.\label{had7b}
\end{eqnarray}
The other two functions, $\Psi_1\left({\bf d}^3;z\right)$ and $\Psi_2\left({
\bf d}^3;z\right)$, can be written by the cyclic permutation of the indices 
$(1,2,3)$ in (\ref{had7b})
\begin{eqnarray}   
\Psi_1\left({\bf d}^3;z\right)=\frac{1}{2\pi}\int_0^{2\pi}\frac{\left(1-
e^{id_2d_3t}\right)dt}{\left(1-e^{id_2t}\right)\left(1-e^{id_3t}\right)\left(
1-z^{d_1}e^{-id_1t}\right)}-1\;,\label{had7c}\\
\Psi_2\left({\bf d}^3;z\right)=\frac{1}{2\pi}\int_0^{2\pi}\frac{\left(1-
e^{id_3d_1t}\right)dt}{\left(1-e^{id_3t}\right)\left(1-e^{id_1t}\right)\left(
1-z^{d_2}e^{-id_2t}\right)}-1\;.\label{had7d}
\end{eqnarray}
\section{The explicit calculation of the entries of the Johnson's matrix}
\label{sect4}
The complexity of expressions (\ref{inter8}) and (\ref{inter10}) in conjunction 
with (\ref{had7b}) -- (\ref{had7d}) makes further evaluation of the diagonal 
elements $a_{kk}\left({\bf d}^3\right)$ of the Johnson's matrix excessively 
difficult. Therefore, we develop in this Section another approach representing 
$a_{kk}\left({\bf d}^3\right)$ as {\em zeroes} of some functions.

Combining (\ref{had7w}) and (\ref{had7b}) we have
\begin{eqnarray}
\frac{1}{2\pi}\int_0^{2\pi}\frac{\left(1-
e^{id_1d_2t}\right)dt}{\left(1-e^{id_1t}\right)\left(1-e^{id_2t}\right)\left(
1-z^{d_3}e^{-id_3t}\right)}-1=\sum_{j=1}^{\infty}z^{b_jd_3}\;,\label{diag1}
\end{eqnarray}
where $b_1,b_2,\ldots$ present all positive integers such that in accordance 
with (\ref{had7v}) every integer $b_jd_3$ is representable by $d_1$ and $d_2$. 

\begin{figure}[h]%[t]
\centerline{\psfig{figure=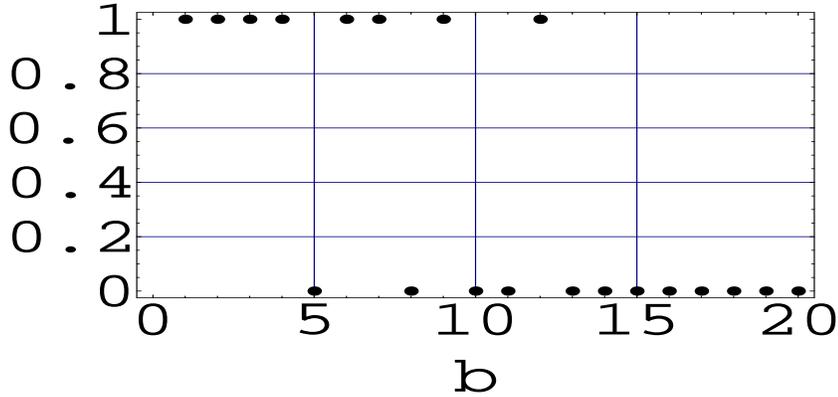,height=5.5cm,width=13cm}}
%width=4.8in}}
\caption{Typical plot of the function $\Xi_3\left(\{23,29,44\};b\right)$ and 
the distribution of its zeroes $b_k$. The diagonal element $a_{33}$ of the 
Johnson's matrix reads $a_{33}\left(23,29,44\right)=5$.}
\label{repr3}
\end{figure}

Fix the index $j=k$ and differentiate with respect to $z^{d_3}$ both sides of 
the equality (\ref{diag1}) $b_k$ times, and take their values at $z=0$
\begin{eqnarray}
\frac{b_k!}{2\pi}\int_0^{2\pi}\frac{\left(1-e^{id_1d_2t}\right)e^{-ib_kd_3t}}
{\left(1-e^{id_1t}\right)\left(1-e^{id_2t}\right)}dt=b_k!\label{diag2}
\end{eqnarray}
Define a new function $\Xi_3\left({\bf d}^3;b\right)$ by
\begin{eqnarray}
\Xi_3\left({\bf d}^3;b\right)=1-\frac{1}{2\pi}\int_0^{2\pi}\frac{\left(1-e^{id_1
d_2t}\right)e^{-ibd_3t}}{\left(1-e^{id_1t}\right)\left(1-e^{id_2t}\right)}dt\;.
\label{diag3}
\end{eqnarray}
Thus $b_k$ is its {\em zero}, $\Xi_3\left({\bf d}^3;b_k\right)=0$, and $a_{33}
\left({\bf d}^3\right)$ is its {\em minimal integer zero} according to 
(\ref{had777})
\begin{eqnarray}
a_{33}\left({\bf d}^3\right)=\min\left\{b_k\;\bracevert\;\Xi_3\left({\bf 
d}^3;b_k\right)=0,\;b_k\in {\mathbb N}\right\}\;.\label{diag4}
\end{eqnarray}
Notice that $2\leq a_{33}\left({\bf d}^3\right)\leq d_1-1$ according to 
\cite{fel04}, that bound the range of $b$ where the first zero $b_1$ does 
appear. In Figure \ref{repr3} we present the typical plot of the function 
$\Xi_3\left({\bf d}^3;b\right)$ for the triple $d_1=23$, $d_2=29$ and $d_3=44$ 
which was considered numerically in \cite{barv02}

The other two diagonal elements, $a_{11}\left({\bf d}^3\right)$ and $a_{22}
\left({\bf d}^3\right)$, can be obtained by the cyclic permutation of the 
indices $(1,2,3)$ in (\ref{diag3}) and 
(\ref{diag4})
\begin{eqnarray}
a_{11}\left({\bf d}^3\right)=\min\left\{b_k\;\bracevert\;\Xi_1\left({\bf 
d}^3;b_k\right)=0,\;b_k\in {\mathbb N}\right\},\;\;
a_{22}\left({\bf d}^3\right)=\min\left\{b_k\;\bracevert\;\Xi_2\left({\bf 
d}^3;b_k\right)=0,\;b_k\in{\mathbb N}\right\}\label{diag5}
\end{eqnarray}
where
\begin{eqnarray}
\Xi_1\left({\bf d}^3;b\right)=1-\frac{1}{2\pi}\int_0^{2\pi}\frac{\left(1-e^{
id_2d_3t}\right)e^{-ibd_1t}}{\left(1-e^{id_2t}\right)\left(1-e^{id_3t}\right)}
dt,\;\;\;\Xi_2\left({\bf d}^3;b\right)=1-\frac{1}{2\pi}\int_0^{2\pi}\frac{
\left(1-e^{id_3d_1t}\right)e^{-ibd_2t}}{\left(1-e^{id_3t}\right)\left(1-e^{i
d_1t}\right)}dt\;.\nonumber
\end{eqnarray}
and $2\leq a_{11}\left({\bf d}^3\right)\leq d_2-1$ , $\;2\leq a_{22}\left({\bf 
d}^3\right)\leq d_1-1$.
 
Performing the calculation for the triple $d_1=23$, $d_2=29$ and $d_3=44$ 
we can find in accordance with (\ref{diag5})
\begin{eqnarray}
a_{11}\left(23,29,44\right)=a_{22}\left(23,29,44\right)=7\;.\nonumber
\end{eqnarray}
Together with $a_{33}\left(23,29,44\right)=5$ these are in full agreement with 
the Johnson's matrix of minimal relations which was found in \cite{fel04}, 
Example 2.
\subsection{The off--digonal elements of the Johnson's matrix}
\label{sect41}
The uniqueness of the Johnson's matrix $((a_{ij}))$ of minimal relations makes 
us possible to determine also its off--digonal elements for non--symmetric 
semigroup ${\sf S}\left({\bf d}^3\right)$. As was shown by Johnson 
\cite{john60}, six off--digonal elements of the matrix $((a_{ij}))$ are related 
by six identities
\begin{eqnarray}
&&a_{21}+a_{31}=a_{11}\;,\;\;\;a_{12}+a_{32}=a_{22}\;,\;\;\;a_{13}+a_{23}=
a_{33}\;,\label{herznon3aaa}\\
&&a_{23}a_{32}=a_{22}a_{33}-d_1\;,\;\;\;a_{13}a_{31}=a_{11}a_{33}-d_2\;,\;\;\;
a_{12}a_{21}=a_{11}a_{22}-d_3\;.\nonumber
\end{eqnarray}
These identities give rise to six quadratic equations
\begin{eqnarray}
d_3a_{23}^2-\left(\langle{\bf a},{\bf d}\rangle-2a_{11}d_1\right)a_{23}+
\left(a_{22}a_{33}-d_1\right)d_2=
d_3a_{13}^2-\left(\langle{\bf a},{\bf d}\rangle-2a_{22}d_2\right)a_{13}+
\left(a_{11}a_{33}-d_2\right)d_1=0\;,\nonumber\\
d_2a_{32}^2-\left(\langle{\bf a},{\bf d}\rangle-2a_{11}d_1\right)a_{32}+
\left(a_{33}a_{22}-d_1\right)d_3=
d_2a_{12}^2-\left(\langle{\bf a},{\bf d}\rangle-2a_{33}d_3\right)a_{12}+
\left(a_{11}a_{22}-d_3\right)d_1=0\;,\nonumber\\
d_1a_{31}^2-\left(\langle{\bf a},{\bf d}\rangle-2a_{22}d_2\right)a_{31}+
\left(a_{33}a_{11}-d_2\right)d_3=
d_1a_{21}^2-\left(\langle{\bf a},{\bf d}\rangle-2a_{33}d_3\right)a_{21}+
\left(a_{22}a_{11}-d_3\right)d_2=0\;,\nonumber
\end{eqnarray}
where $\langle{\bf a},{\bf d}\rangle$ is already defined in (\ref{frob3}).
Notice that all these equations have common discriminant
$$
\langle{\bf a},{\bf d}\rangle^2-4\left(a_{11}a_{22}d_1d_2+a_{22}a_{33}d_2d_3+
a_{33}a_{11}d_3d_1\right)+4d_1d_2d_3\;,
$$
which can be recognized as $J^2\left({\bf d}^3\right)$ defined in (\ref{frob3}).
As was shown in \cite{fel04} its square root, $J\left({\bf d}^3\right)$, is a 
positive integer,  
\begin{eqnarray}
J\left({\bf d}^3\right)=\arrowvert a_{12}a_{23}a_{31}-a_{13}a_{32}a_{21}
\arrowvert\geq 1\;.\label{pos1}
\end{eqnarray}
Therefore the solutions of all six quadratic equations are always rational 
numbers. 
\begin{eqnarray}
&&a_{23}^{\pm}=\frac{1}{2d_3}\left(\langle{\bf a},{\bf d}\rangle\pm J\left({\bf 
d}^3\right)-2a_{11}d_1\right)\;,\;\;\;a_{32}^{\pm}=\frac{1}{2d_2}\left(\langle{
\bf a},{\bf d}\rangle\pm J\left({\bf d}^3\right)-2a_{11}d_1\right)\;,\nonumber\\
&&a_{31}^{\pm}=\frac{1}{2d_1}\left(\langle{\bf a},{\bf d}\rangle\pm J\left({\bf 
d}^3\right)-2a_{22}d_2\right)\;,\;\;\;a_{13}^{\pm}=\frac{1}{2d_3}\left(\langle{
\bf a},{\bf d}\rangle\pm J\left({\bf d}^3\right)-2a_{22}d_2\right)\;,\nonumber\\
&&a_{12}^{\pm}=\frac{1}{2d_2}\left(\langle{\bf a},{\bf d}\rangle\pm J\left({\bf
d}^3\right)-2a_{33}d_3\right)\;,\;\;\;a_{21}^{\pm}=\frac{1}{2d_1}\left(\langle{
\bf a},{\bf d}\rangle\pm J\left({\bf d}^3\right)-2a_{33}d_3\right)\;.
\label{pos2}
\end{eqnarray}
Show that one of two rational roots associated with every quadratic equation is 
always positive integer while another is necessarily not integer. Making use 
of the identity \cite{fel04}, formula (135),
\begin{eqnarray}
\frac{1}{2}\left[\langle{\bf a},{\bf d}\rangle\pm J\left({\bf d}^3\right)\right]
=a_{11}a_{22}a_{33}+\frac{1}{2}\left(a_{12}a_{23}a_{31}+a_{13}a_{32}a_{21}\pm
\arrowvert a_{12}a_{23}a_{31}-a_{13}a_{32}a_{21}\arrowvert\right)\;,\nonumber
\end{eqnarray}
and the identities (\ref{herznon3aaa}), we have
\begin{eqnarray}
a_{23}^{+}=a_{23}\;,\;\;\;a_{32}^{+}=a_{32}\;,\;\;\;a_{31}^{+}=a_{31}\;,\;\;\;
a_{13}^{+}=a_{13}\;,\;\;\;a_{12}^{+}=a_{12}\;,\;\;\;a_{21}^{+}=a_{21}\;,\;\;
\label{pos3}\\
a_{23}^{-}=a_{32}\frac{d_2}{d_3}\;,\;\;\;a_{32}^{-}=a_{23}\frac{d_3}{d_2}\;,
\;\;\;a_{31}^{-}=a_{13}\frac{d_3}{d_1}\;,\;\;\;a_{13}^{-}=a_{31}\frac{d_1}{d_3}
\;,\;\;\;a_{12}^{-}=a_{21}\frac{d_1}{d_2}\;,\;\;\;a_{21}^{-}=a_{12}\frac{d_2}
{d_1}\;.\;\;\label{pos4}
\end{eqnarray}
There is only one way to satisfy the uniqueness of the Johnson's matrix 
$((a_{ij}))$ of minimal relations comprised exclusively of integer entries 
$a_{ij}$ if we require that $a_{ij}^{-}=a_{ij}^{+}=a_{ij}$. However, as one can 
see from (\ref{pos2}), this leads to $J\left({\bf d}^3\right)=0$ that 
contradicts (\ref{pos1}).
\section*{Acknowledgement}
The author thanks A. Juhasz for usefull discussions.

\end{document}